\documentclass[12pt]{amsart}
\usepackage{amssymb}
\usepackage{amsfonts}

\newtheorem{theorem}{Theorem}[section]
\theoremstyle{plain}

\newtheorem{corollary}[theorem]{Corollary}

\newtheorem{example}[theorem]{Example}

\newtheorem{lemma}[theorem]{Lemma}

\newtheorem{proposition}[theorem]{Proposition}
\newtheorem{remark}[theorem]{Remark}

\numberwithin{equation}{section}

\input{tcilatex}

\begin{document}
\title{C*-crossed-products by an order-two automorphism}
\author{Man-Duen Choi}
\address{Department of Mathematics\\
University of Toronto}
\email{choi@math.toronto.edu}
\author{Fr\'{e}d\'{e}ric Latr\'{e}moli\`{e}re}
\email{frederic@math.toronto.edu}
\urladdr{http://www.math.toronto.edu/frederic}
\date{June 2006}

\begin{abstract}
We describe the representation theory of C*-crossed-products of a unital
C*-algebra $A$ by the cyclic group of order 2. We prove that there are two
main types of irreducible representations for the crossed-product: those
whose restriction to $A$ is irreducible and those who are the sum of two
unitarily unequivalent representations of $A$. We characterize each class in
term of the restriction of the representations to the fixed point
C*-subalgebra of $A$. We apply our results to compute the $K$-theory of
several crossed-products of the free group on two generators.
\end{abstract}

\maketitle

\section{Introduction}

This paper explores the structure of the representation theory of
C*-crossed-products \cite{Zeller-Meier68}\ of unital C*-algebras by
order-two automorphisms. We show that irreducible representations of the
C*-crossed-products $A\rtimes \mathbb{Z}_{2}$ of a unital C*-algebra $A$ by $%
\mathbb{Z}_{2}$ fall in two categories: either their restriction to $A$ is
already irreducible, or it is the direct sum of two irreducible
representations of $A$, related together by the automorphism and not
unitarily equivalent to each other.

The paper starts with the given data of a unital C*-algebra $A$ and an
order-two automorphism $\sigma $ of $A$. The C*-crossed-product $A\rtimes
_{\sigma }\mathbb{Z}_{2}$ is the C*-algebra generated by $A$ and a unitary $%
W $ satisfying the following universal property: given any unital *-morphism 
$\psi :A\longrightarrow B$ for some unital C*-algebra $B$ such that $B$
contains a unitary $u$ such that $u^{2}=1$ and $u\psi (a)u^{\ast }=\psi
\circ \sigma (a)$ for all $a\in A$, then $\psi $ extends uniquely to $%
A\rtimes _{\sigma }\mathbb{Z}_{2}$ with $\psi (W)=u$. The general
construction of $A\rtimes _{\sigma }\mathbb{Z}_{2}$ can be found in \cite%
{Zeller-Meier68}. In particular, $W^{2}=1$ (so $W=W^{\ast }$ since $W$ is
unitary) and $WaW^{\ast }=\sigma (a)$ for all $a\in A$. We call the unitary $%
W$ the canonical unitary of $A\rtimes _{\sigma }\mathbb{Z}_{2}$. Proposition
(\ref{Rep1}) in this paper will offer an alternative description of $%
A\rtimes _{\sigma }\mathbb{Z}_{2}$.

The question raised in this paper is: what is the connection between the
representation theory of $A\rtimes _{\sigma }\mathbb{Z}_{2}$ and the
representation theory of $A$? Of central importance is the fixed point
C*-algebra $A_{1}$ for $\sigma $ defined by $A_{1}=\left\{ a\in A:\sigma
(a)=a\right\} $ and the natural decomposition $A=A_{1}+A_{-1}$ where $%
A_{-1}=\left\{ a\in A:\sigma (a)=-a\right\} $, with $A_{1}\cap
A_{-1}=\left\{ 0\right\} $. We obtain a complete description of the
irreducible representations of $A\rtimes _{\sigma }\mathbb{Z}_{2}$ from the
representation theory of $A$ and $A_{1}$.

Note that, if we considered the crossed-product $A\rtimes _{\sigma }\mathbb{Z%
}$ instead of $A\rtimes _{\sigma }\mathbb{Z}_{2}$, then our work applies as
well thanks to a simple observation made at the end of the first section of
this paper.

The rest of the paper focuses on applications to examples. We are interested
in several natural order-two automorphisms of the full C*-algebra of free
group $\mathbb{F}_{2}$, namely the universal C*-algebra generated by two
unitaries $U$ and $V$. We define the automorphism $\alpha $ by $\alpha
(U)=U^{\ast }$ and $\alpha (V)=V^{\ast }$, while $\beta $ is the
automorphism defined by $\beta (U)=-U$ and $\beta (V)=-V$. We compute in
this paper the $K$-theory of the C*-crossed-products for these two
automorphisms, relying in part on our structure theory for their
representations. A third natural automorphism, $\gamma $, is defined
uniquely by $\gamma (U)=V$ and $\gamma (V)=U$. It is the subject of the
companion paper \cite{CLat06} which emphasizes the interesting structure of
the associated fixed point C*-algebra and uses different techniques from the
representation approach of this paper.

\section{Representation theory of the crossed-products}

In this section, we derive several general results on the irreducible
representations of the crossed-product C*-algebra $A\rtimes _{\sigma }%
\mathbb{Z}_{2}$ where $\sigma $ is an order-2 automorphism of the unital
C*-algebra $A$. We recall that $A\rtimes _{\sigma }\mathbb{Z}_{2}$ is the
universal C*-algebra generated by $A$ and a unitary $W$ such that $W^{2}=1$
and $WaW^{\ast }=\sigma (a)$.

\subsection{Representations from the algebra}

A central feature of the crossed-products by finite groups is their
connection with the associated fixed point C*-algebra \cite{Rieffel80}. In
our case, the following easy lemma will prove useful:

\begin{lemma}
Let $A$ be a unital C*-algebra and $\sigma $ an order-2 automorphism of $A$.
The set $A_{1}=\left\{ a+\sigma (a):a\in A\right\} $ is the fixed point
C*-algebra of $A$ for $\sigma $ and the set $A_{-1}=\left\{ a-\sigma
(a):a\in A\right\} $ is the space of elements $b\in A$ such that $\sigma
(b)=-b$. Then $A=A_{1}+A_{-1}$ and $A_{1}\cap A_{-1}=\left\{ 0\right\} $.
\end{lemma}

\begin{proof}
If $a$ is any element in $A$ then $a+\sigma (a)$ (resp. $a-\sigma (a)$) is a
fixed point for $\sigma $ (resp. an element $b\in A$ such that $\sigma
(b)=-b $). Conversely, let $x\in A$: then $x=\frac{1}{2}\left( x+\sigma
(x)\right) +\frac{1}{2}\left( x-\sigma (x)\right) $. If $x$ is $\sigma $%
-invariant then $x-\sigma (x)=0$ so $x=\frac{1}{2}\left( x+\sigma (x)\right) 
$ indeed, and thus the fixed point C*-algebra is $A_{1}$ (and similarly $%
\left\{ b\in A:\sigma (b)=-b\right\} =A_{-1}$).\ Of course, if $a\in
A_{1}\cap A_{-1}$ then $\sigma (a)=a=-a$ so $a=0$.
\end{proof}

We exhibit a simple algebraic description of the crossed-product:

\begin{proposition}
\label{Rep1}Let $\sigma $ be an order 2-automorphism of a unital C*-algebra $%
A$. Then the C*-crossed-product $A\rtimes _{\sigma }\mathbb{Z}_{2}$ is
*-isomorphic to:%
\begin{equation*}
\left\{ \left[ 
\begin{array}{cc}
a & b \\ 
\sigma (b) & \sigma (a)%
\end{array}%
\right] :a,b\in A\right\} \subseteq M_{2}(A)
\end{equation*}%
via the following isomorphism: $a\in A\mapsto \left[ 
\begin{array}{cc}
a & 0 \\ 
0 & \sigma (a)%
\end{array}%
\right] $ and $W\mapsto \left[ 
\begin{array}{cc}
0 & 1 \\ 
1 & 0%
\end{array}%
\right] $ where $W$ is the canonical unitary of $A\rtimes _{\sigma }\mathbb{Z%
}_{2}$.
\end{proposition}

\begin{proof}
Let $\psi :a\in A\mapsto \left[ 
\begin{array}{cc}
a & 0 \\ 
0 & \sigma (a)%
\end{array}%
\right] \in M_{2}(A)$ and set $\psi (W)=\left[ 
\begin{array}{cc}
0 & 1 \\ 
1 & 0%
\end{array}%
\right] \in M_{2}(A)$. Since $\psi (W)\psi (a)\psi (W)=\psi (\sigma (a))$ we
deduce by universality that $\psi $ extends to a (unique)\ *-automorphism of 
$A\rtimes _{\sigma }\mathbb{Z}_{2}$ valued in $M_{2}(A)$. Now, let $c\in
A\rtimes _{\sigma }\mathbb{Z}_{2}$. By construction of $A\rtimes _{\sigma }%
\mathbb{Z}_{2}$, there exists a sequence $\left( a_{n}+b_{n}W\right) _{n\in 
\mathbb{N}}$ with $a_{n},b_{n}\in A$ such that $c=\lim_{n\rightarrow \infty
}a_{n}+b_{n}W$ in $A\rtimes _{\sigma }\mathbb{Z}_{2}$. Now, $\psi
(a_{n}+b_{n}W)=\left[ 
\begin{array}{cc}
a_{n} & b_{n} \\ 
\sigma (b_{n}) & \sigma (a_{n})%
\end{array}%
\right] $ for all $n\in \mathbb{N}$, and converges to $\psi (c)=\left[ 
\begin{array}{cc}
c_{11} & c_{12} \\ 
c_{21} & c_{22}%
\end{array}%
\right] $ when $n\rightarrow \infty $. In particular, $(a_{n})_{n\in \mathbb{%
N}}$ converges to $c_{11}\in A$ and $(b_{n})_{n\in \mathbb{N}}$ converges to 
$c_{12}\in A$. Consequently, $c=c_{11}+c_{12}W$. Hence $A+AW$ is a closed
dense *-subalgebra of $A\rtimes _{\sigma }\mathbb{Z}_{2}$ and thus $A\rtimes
_{\sigma }\mathbb{Z}_{2}=A+AW$.

Moreover, if $\psi (c)=0$ then, writing $c=a+bW$, by definition of $\psi $,
we get $\psi (c)=\left[ 
\begin{array}{cc}
a & b \\ 
\sigma (b) & \sigma (a)%
\end{array}%
\right] =0$ so $a=b=0$ hence $c=0$. Thus $\psi $ is a *-isomorphism from $%
A\rtimes _{\sigma }\mathbb{Z}_{2}$ onto the C*-algebra $\left\{ \left[ 
\begin{array}{cc}
a & b \\ 
\sigma (b) & \sigma (a)%
\end{array}%
\right] :a,b\in A\right\} \subseteq M_{2}(A)$.

This concludes our proof.
\end{proof}

In other words, the abstract canonical unitary $W$ of $A\rtimes _{\sigma }%
\mathbb{Z}_{2}\ $can be replaced by the concrete unitary $\left[ 
\begin{array}{cc}
0 & 1 \\ 
1 & 0%
\end{array}%
\right] $ and $A\rtimes _{\sigma }\mathbb{Z}_{2}$ can be seen as the
C*-algebra $\psi (A)+\psi (A)\left[ 
\begin{array}{cc}
0 & 1 \\ 
1 & 0%
\end{array}%
\right] $ in $M_{2}(A)$ with $\psi :a\in A\mapsto \left[ 
\begin{array}{cc}
a & 0 \\ 
0 & \sigma (a)%
\end{array}%
\right] $. Equivalently, the *-subalgebra $A+AW$ in $A\rtimes _{\sigma }%
\mathbb{Z}_{2}$ is in fact equal to $A\rtimes _{\sigma }\mathbb{Z}_{2}$.

From the algebraic description of Proposition (\ref{Rep1}) we get a family
of representations of the crossed-product described in the following
proposition. These representations are in fact induced representations from
the sub-C*-algebra $A$ to the C*-algebra $A\rtimes _{\sigma }\mathbb{Z}_{2}$
in the sense of \cite{Rieffel74}.

\begin{proposition}
\label{Rep0}Let $A$ be a unital C*-algebra and $\sigma $ be an order two
automorphism of $A$. Let $W$ be the canonical unitary of the crossed-product 
$A\rtimes _{\sigma }\mathbb{Z}_{2}$ such that $WaW=\sigma (a)$ for all $a\in
A$. Then, for each representation $\pi $ of $A$ on some Hilbert space $%
\mathcal{H}$ there exists a representation $\widetilde{\pi }$ of $A\rtimes
_{\sigma }\mathbb{Z}_{2}$ on $\mathcal{H}\oplus \mathcal{H}$\ defined by $%
\widetilde{\pi }(a)=\pi (a)\oplus \pi \circ \sigma (a)$ for all $a\in A$ and 
$\widetilde{\pi }(W)=\left[ 
\begin{array}{cc}
0 & 1 \\ 
1 & 0%
\end{array}%
\right] $. The following are moreover equivalent:

\begin{itemize}
\item The representation $\widetilde{\pi }$ is irreducible,

\item The representation $\pi $ is irreducible and not unitarily equivalent
to $\pi \circ \sigma $,

\item There does not exist a unitary $U\in \mathcal{B}(\mathcal{H})$ such
that $U\pi U^{\ast }=\pi \circ \sigma $ and $U^{2}=1$.
\end{itemize}

If $\pi $ is a faithful representation of $A$ then $\widetilde{\pi }$ is
faithful for $A\rtimes _{\sigma }\mathbb{Z}_{2}$. In particular, if $A$ has
a faithful representation which is a direct sum of finite representations,
so does $A\rtimes _{\sigma }\mathbb{Z}_{2}$.
\end{proposition}

\begin{proof}
Let $\pi $ be a given representation of $A$. Then by setting $\widetilde{\pi 
}(a)=\pi (a)\oplus \pi (\sigma (a))$ and $\widetilde{\pi }(W)=\left[ 
\begin{array}{cc}
0 & 1 \\ 
1 & 0%
\end{array}%
\right] $, we define a *-representation of $A\rtimes _{\sigma }\mathbb{Z}%
_{2} $ by universality of $A\rtimes _{\sigma }\mathbb{Z}_{2}$. In fact, $%
\widetilde{\pi }=\left[ 
\begin{array}{cc}
\pi & 0 \\ 
0 & \pi%
\end{array}%
\right] \circ \psi $ where $\psi $ is the isomorphism of Proposition (\ref%
{Rep1}).

Let us now assume that $\pi $ is irreducible and not unitarily equivalent to 
$\pi \circ \sigma $. Assume $V$ is an operator commuting with $\widetilde{%
\pi }$. Then, since $V$ commutes with $\widetilde{\pi }(W)$, we have $V=%
\left[ 
\begin{array}{cc}
a & b \\ 
b & a%
\end{array}%
\right] $ for some $a,b\in \mathcal{B}\left( \mathcal{H}\right) $. Now,
since $V$ commutes with $\pi \oplus \left( \pi \circ \sigma \right) $ we
conclude that $a$ commutes with $\pi $ and, as $\pi $ is irreducible, this
implies that $V=\left[ 
\begin{array}{cc}
\lambda 1 & b \\ 
b & \lambda 1%
\end{array}%
\right] $ for some $\lambda \in \mathbb{C}$ and where $1$ is the identity on 
$\mathcal{H}$. Hence, $V$ commutes with $\widetilde{\pi }$ if and only if $%
\left[ 
\begin{array}{cc}
0 & b \\ 
b & 0%
\end{array}%
\right] $ does. Now, if $\left[ 
\begin{array}{cc}
0 & b \\ 
b & 0%
\end{array}%
\right] $ commutes with $\widetilde{\pi }(A\rtimes _{\sigma }\mathbb{Z}_{2})$
then so does its square $\left[ 
\begin{array}{cc}
b^{2} & 0 \\ 
0 & b^{2}%
\end{array}%
\right] $. Hence again $b^{2}=\mu 1$ by irreducibility of $\pi $, and up to
replacing $b$ by $\frac{1}{2}\left( b+b^{\ast }\right) $ we can assume that $%
b$ is self-adjoint and thus $\mu \geq 0$.

Assume that $\mu \not=0$. Set $u=\left( \sqrt[2]{\mu }\right) ^{-1}b$: then $%
u=u^{\ast }$ and $u^{2}=1$ so $u$ is a unitary. Moreover, as $\frac{1}{\sqrt[%
2]{\mu }}V=\left[ 
\begin{array}{cc}
0 & u \\ 
u & 0%
\end{array}%
\right] $ commutes with $\pi (c)\oplus \left( \pi \circ \sigma (c)\right) $
for all $c\in A$, we check that $u\pi (\sigma (c))=$ $\pi (c)u$ so $u^{\ast
}\pi (c)u=\pi (\sigma (c))$ for all $c\in A$. Hence, we have reached a
contradiction as we assumed that $\pi $ is not unitarily equivalent to $\pi
\circ \sigma $. Therefore $\mu =0$ and thus $V=\lambda \left( 1\oplus
1\right) $, so $\widetilde{\pi }$ is irreducible.

Conversely, if there exists a unitary $u$ such that $u^{2}=1$ and $u\pi
u^{\ast }=\pi \circ \sigma $, then the operator $V=\left[ 
\begin{array}{cc}
0 & u \\ 
u & 0%
\end{array}%
\right] $ commutes with $\widetilde{\pi }$ so $\widetilde{\pi }$ is not
irreducible.

On the other hand, if $\pi $ is reducible, then let $p$ be a nontrivial
projection of $\mathcal{H}$ such that $p\pi =\pi p$. Then $p\oplus p$ is a
nontrivial projection commuting with $\widetilde{\pi }$ as can easily been
checked (it is obvious on $\widetilde{\pi }(A)$ and easy for $\widetilde{\pi 
}(W)$). Hence $\widetilde{\pi }$ is reducible as well. This proves the first
two equivalence.

Now, we observe that $\pi $ is unitarily equivalent to $\pi \circ \sigma $
if and only if there exists a unitary $u$ with $u^{2}=1$ such that $u\pi
u^{\ast }=\pi \circ \sigma $. One implication is trivial; let us check the
easy other one. Let $v$ be a unitary such that $v\pi v^{\ast }=\pi \circ
\sigma $. Then $v^{2}\pi v^{\ast 2}=\pi \circ \sigma ^{2}=\pi $. Hence, as $%
\pi $ is irreducible, $v^{2}=\exp \left( 2i\pi \theta \right) 1$ for some $%
\theta \in \lbrack 0,1)$. Hence, $u=\exp \left( -i\pi \theta \right) v$
satisfies both $u^{2}=1$ and $\pi \circ \sigma =u\pi u^{\ast }$.
\end{proof}

Proposition (\ref{Rep0}) describes a family of representations and gives us
a criterion for their irreducibility. Conversely, given an irreducible
representation of $A\rtimes _{\sigma }\mathbb{Z}_{2}$, what can be said
about its structure relative to the representation theory of $A$ and its
fixed point algebra $A_{1}$?\ This is the matter of the next section, which
establishes a sort of converse for Proposition (\ref{Rep0}).

\subsection{Irreducible Representations}

We will use the following lemma:

\begin{lemma}
\label{central}Let $\mathcal{H}$ be a Hilbert space. Let $A,B$ be two
bounded linear operators on $\mathcal{H}$ such that $BTA=ATB$ for all
bounded linear operators $T$ on $\mathcal{H}$. Then $A$ and $B$ are linearly
dependant.
\end{lemma}

\begin{proof}
The result is obvious if $A=0$ or $B=0$, so we assume henceforth that $%
A\not=0$ and $B\not=0$. Let $\gamma \in \mathcal{H}$ such that $A\gamma
\not=0$. Assume that there exists $x_{0}\in \mathcal{H}$ such that $\left\{
Ax_{0},Bx_{0}\right\} $ is linearly independent. Then let $T$ be any bounded
linear operator such that $T(Ax_{0})=0$ and $T(Bx_{0})=\gamma $. Such a $T$
is well-defined by the Hahn-Banach theorem. But then $0=BTAx_{0}=ATBx_{0}=A%
\gamma $ which is a contradiction.\ Hence for all $x\in \mathcal{H}$ there
exists $\lambda _{x}\in \mathbb{C}$ such that $Bx=\lambda _{x}Ax$.

Now, let $y\in \mathcal{H}$. Let $T$ be any bounded operator on $\mathcal{H}$
such that $TA\gamma =y$. Then we compute:%
\begin{equation*}
By=BTA\gamma =ATB\gamma =AT(\lambda _{\gamma }A\gamma )=\lambda _{\gamma }Ay%
\text{.}
\end{equation*}%
Hence $B=\lambda _{\gamma }A$. This concludes our theorem.
\end{proof}

Note that we can prove similarly:

\begin{lemma}
\label{central2}Let $A$, $B$ be two bounded operators on a Hilbert space $%
\mathcal{H}$ and assume that for all bounded operators $T$ of $\mathcal{H}$
we have:%
\begin{equation*}
ATA^{\ast }=BTB^{\ast }\text{.}
\end{equation*}%
Then there exists $\theta \in \lbrack 0,1)$ such that $B=\exp \left( 2i\pi
\theta \right) A$.
\end{lemma}

Either lemma can be used to prove the following description of the structure
of irreducible representations of $A\rtimes _{\sigma }\mathbb{Z}_{2}$. This
theorem is the main result of this paper, and shows that any irreducible
representation of $A\rtimes _{\sigma }\mathbb{Z}_{2}$ is build from either a
single unitary representation of $A$ (and is then just an extension of it)
or from two non-equivalent irreducible representations of $A$.

\begin{theorem}
\label{Rep}Let $\sigma $ be an order--two--automorphism of a unital
C*-algebra $A$. We denote by $W$ the canonical unitary of the
C*-crossed-product $A\rtimes _{\sigma }\mathbb{Z}_{2}$ such that $WaW=\sigma
(a)$ for all $a\in A$.

Let $\pi $ be an irreducible representation of $A\rtimes _{\sigma }\mathbb{Z}%
_{2}$ on a Hilbert space $\mathcal{H}$. Let $\pi ^{\prime }$ be the
restriction of $\pi $ to $A$ and $\pi ^{\prime \prime }$ be the restriction
of $\pi $ to the fixed point C*-algebra $A_{1}$. Then one and only one of
the following two alternatives hold:

\begin{enumerate}
\item the operator $\pi (W)$ is either the identity $\limfunc{Id}$ or $-%
\limfunc{Id}$ and $\pi (A\rtimes _{\sigma }\mathbb{Z}_{2})=\pi ^{\prime
}(A)=\pi ^{\prime \prime }(A_{1})$.

\item the spectrum of $\pi (W)$ is $\left\{ -1,1\right\} $. Then $\mathcal{H}%
=\mathcal{H}_{1}\oplus \mathcal{H}_{-1}$ where $\mathcal{H}_{\varepsilon }$
is the spectral Hilbert space of $\pi (W)$ for the eigenvalue $\varepsilon $%
. With this decomposition of $\mathcal{H}$, we have $\pi (W)=\left[ 
\begin{array}{cc}
1 & 0 \\ 
0 & -1%
\end{array}%
\right] $. Let us write $\pi ^{\prime }(a)=\left[ 
\begin{array}{cc}
\alpha (a) & \beta (a) \\ 
\gamma (a) & \delta (a)%
\end{array}%
\right] $ for $a\in A$. Then $\alpha ,\delta $ restrict to irreducible
representations of $A_{1}$, and $\alpha (A_{-1})=\delta (A_{-1})=\left\{
0\right\} $. Moreover, $\beta (A_{1})=\gamma (A_{1})=\left\{ 0\right\} $.

Furthermore, the representation $\pi ^{\prime }$ is irreducible if and only
if $\alpha $ and $\delta $ are not unitarily equivalent.
\end{enumerate}
\end{theorem}

\begin{proof}
Let $\pi $ be an irreducible representation of $A\rtimes _{\sigma }\mathbb{Z}%
_{2}$ on $\mathcal{H}$. Let $w=\pi (W)$. Since $w$ is unitary and $w^{2}=1$,
the spectrum of $w$ is either $\left\{ -1,1\right\} $ or $w=1$ or $w=-1$. In
the latter two cases, $w$ commutes with $\pi \left( A\rtimes _{\sigma }%
\mathbb{Z}_{2}\right) $. Since $A\rtimes _{\sigma }\mathbb{Z}_{2}=A+AW$ from
Proposition (\ref{Rep1}) we have $\pi \left( A\rtimes _{\sigma }\mathbb{Z}%
_{2}\right) =\pi ^{\prime }(A)+\pi ^{\prime }(A)w=\pi ^{\prime }(A)$ (as $%
w=\pm 1$). Thus as $\pi $ is irreducible, so is $\pi ^{\prime }$. Moreover,
since $w\pi ^{\prime }(a)w=\pi ^{\prime }(a)=\pi ^{\prime }\circ \sigma (a)$%
, we see that $\pi ^{\prime }$ is null on $A_{-1}$ and thus $\pi ^{\prime
}=\pi ^{\prime \prime }$. Conversely if $\pi (A_{-1})=0$ then $w$ must
commute with $\pi (A)=\pi (A_{1})$ and thus with $\pi (A\rtimes _{\sigma }%
\mathbb{Z}_{2})=\pi (A)+\pi (A)w$. Therefore, as $\pi $ is irreducible, $w$
is scalar, and as $w$ unitary and $w^{2}=1$ we conclude $w$ is $1$ or $-1$.

Assume now that the unitary $w$ has spectrum $\left\{ -1,1\right\} $. Write $%
\mathcal{H}=\mathcal{H}_{1}\oplus \mathcal{H}_{-1}$ accordingly. In this
decomposition, we have 
\begin{equation*}
w=\left[ 
\begin{array}{cc}
1 & 0 \\ 
0 & -1%
\end{array}%
\right] \text{ and }\pi (a)=\left[ 
\begin{array}{cc}
\alpha (a) & \beta (a) \\ 
\gamma (a) & \delta (a)%
\end{array}%
\right]
\end{equation*}%
where $\alpha ,\beta ,\gamma ,\delta $ are linear maps on $A$. Thus:%
\begin{equation*}
w\pi (a)w^{\ast }=\left[ 
\begin{array}{cc}
\alpha (a) & -\beta (a) \\ 
-\gamma (a) & \delta (a)%
\end{array}%
\right] \text{.}
\end{equation*}
In particular, if $a\in A_{-1}$ then $\pi \circ \sigma (a)=-\pi (a)$ so $%
\alpha (a)=-\alpha (a)=0$. Since $A=A_{1}\oplus A_{-1}$ as a vector space,
we conclude that $\alpha (a)\in \alpha (A_{1})$ for all $a\in A$. Similarly $%
\delta (a)\in \delta (A_{1})$, $\beta (a)\in \beta (A_{-1})$ and $\gamma
(a)\in \gamma (A_{-1})$ for all $a\in A$ while $\gamma (A_{1})=\beta
(A_{1})=\left\{ 0\right\} $.

Consequently, $\pi ^{\prime \prime }=\alpha \oplus \beta $ and $\alpha
,\beta $ are representations of $A_{1}$ (but not of $A$).

We observe that $A\rtimes _{\sigma }\mathbb{Z}_{2}=A+AW$ by Proposition (\ref%
{Rep1}), so: 
\begin{equation*}
\pi \left( A\rtimes _{\sigma }\mathbb{Z}_{2}\right) =\left\{ \left[ 
\begin{array}{cc}
\alpha (a_{1})+\alpha (a_{2}) & \beta (a_{1})-\beta (a_{2}) \\ 
\gamma (a_{1})+\gamma (a_{2}) & \delta (a_{1})-\delta (a_{2})%
\end{array}%
\right] :a_{1},a_{2}\in A\right\}
\end{equation*}%
(note that $w$ is given in this form by $a_{1}=1$ and $a_{2}=0$, since $1\in
A_{1}$ so $\beta (a_{1})=\gamma (a_{1})=0$). Now, $\alpha (a)\in \alpha
(A_{1})$ for all $a\in A$, so $\left\{ \alpha (a_{1})+\alpha \left(
a_{2}\right) :a_{1},a_{2}\in A\right\} $ is the set $\alpha (A_{1})$.
Furthermore, since $\pi $ is irreducible, we have $\pi \left( A\rtimes
_{\sigma }\mathbb{Z}_{2}\right) ^{\prime \prime }=\mathcal{B}\left( \mathcal{%
H}\right) $, i.e. the range of $\pi $ is $\limfunc{SOT}$-dense, and in
particular $\alpha (A_{1})$ is $\limfunc{SOT}$-dense in $\mathcal{B}\left( 
\mathcal{H}_{1}\right) $, so $\alpha $ is an irreducible representation of $%
A_{1}$ on $\mathcal{H}_{1}$. The same applies to $\delta $.

We now distinguish according to the two following cases: either $\alpha $
and $\delta $ are unitarily equivalent as representations of $A_{1}$ or they
are not.

Assume that $\alpha $ and $\delta $ are not unitarily equivalent. Let us
assume $P$ is a projection which commutes with $\pi ^{\prime }$. Then in
particular, $P$ commutes with $\pi ^{\prime \prime }$. Writing $P=\left[ 
\begin{array}{cc}
p_{11} & p_{12} \\ 
p_{21} & p_{22}%
\end{array}%
\right] $, this gives the relations:%
\begin{eqnarray*}
\left[ 
\begin{array}{cc}
\alpha & 0 \\ 
0 & \delta%
\end{array}%
\right] \left[ 
\begin{array}{cc}
p_{11} & p_{12} \\ 
p_{21} & p_{22}%
\end{array}%
\right] &=&\left[ 
\begin{array}{cc}
\alpha p_{11} & \alpha p_{12} \\ 
\delta p_{21} & \delta p_{22}%
\end{array}%
\right] \\
\left[ 
\begin{array}{cc}
p_{11} & p_{12} \\ 
p_{21} & p_{22}%
\end{array}%
\right] \left[ 
\begin{array}{cc}
\alpha & 0 \\ 
0 & \delta%
\end{array}%
\right] &=&\left[ 
\begin{array}{cc}
p_{11}\alpha & p_{12}\delta \\ 
p_{21}\alpha & p_{22}\delta%
\end{array}%
\right]
\end{eqnarray*}%
Hence, since both $\alpha $ and $\delta $ are irreducible, we deduce that $%
p_{11},p_{22}$ are scalar. Now, as $P$ is a projection, $p_{11}=p_{11}^{\ast
}$ and $\left\Vert p_{11}\right\Vert \leq 1$ so $p_{11}\in \lbrack -1,1]$.
Again since $P=P^{\ast }=P^{2}$ we have $p_{12}=p_{21}^{\ast }$ and $%
p_{12}p_{12}^{\ast }+p_{11}^{2}=p_{11}\in \mathbb{R}$. Assume $\lambda
=p_{11}\left( 1-p_{11}\right) \not=0$. Since $p_{11}\in \lbrack -1,1]$ we
have $\lambda \in \lbrack 0,1]$. Then $\nu =\frac{1}{\sqrt[2]{\lambda }}%
p_{12}$ is a unitary operator and since $p_{12}^{\ast }\alpha p_{12}=\delta $%
, we obtain $\nu \alpha \nu ^{\ast }=\delta $. This contradicts our
assumption that $\alpha $ and $\delta $ are not unitarily equivalent. Hence $%
\lambda =0$ and so $p_{11}=1$ or $0$ and $p_{12}=0$ (since $%
p_{12}p_{12}^{\ast }=0$). Now, again since $P$ is a projection, $%
p_{22}^{2}+p_{12}p_{12}^{\ast }=p_{22}$ yet $p_{12}=0$ and $p_{22}$ is a
scalar so $p_{22}=0$ or $1$ as well. Thus, in the decomposition $\mathcal{H}=%
\mathcal{H}_{1}\oplus \mathcal{H}_{-1}$ the projection $P$ is either $%
0\oplus 0$, $1\oplus 0$, $0\oplus 1$ or $1\oplus 1$.

Now, the first part of this proof established that $\pi (W)$ must be scalar
if $\pi $ is irreducible and $\pi (A)=\pi (A_{1})$. Since we assume that $%
\pi (W)$ is not scalar, we conclude that $\pi (A)\not=\pi (A_{1})$.
Consequently, there exists $a_{0}\in A\backslash A_{1}$ such that $\pi
(a_{0})$ is not diagonal in the decomposition $\mathcal{H}_{1}\oplus 
\mathcal{H}_{-1}$. Thus $\pi (a_{0})$ does not commute with $\left[ 
\begin{array}{cc}
1 & 0 \\ 
0 & 0%
\end{array}%
\right] $ and $\left[ 
\begin{array}{cc}
0 & 0 \\ 
0 & 1%
\end{array}%
\right] $. So $P$ is scalar, and thus $\pi ^{\prime }$ is irreducible. Note
that $\pi (W)\pi ^{\prime }(a)\pi (W)=\pi ^{\prime }\circ \sigma (a)$ for
all $a\in A$, so $\pi ^{\prime }$ is unitarily equivalent to $\pi ^{\prime
}\circ \sigma $.

Conversely, assume $\alpha $ and $\delta $ are unitarily equivalent. Thus,
there exists a unitary $u\in \mathcal{B}(\mathcal{H}_{1},\mathcal{H}_{-1})$
such that $\alpha =u\beta u^{\ast }$. By conjugating $\pi $ with $u^{\prime
}=\left[ 
\begin{array}{cc}
1 & 0 \\ 
0 & u%
\end{array}%
\right] $, we obtain $u^{\prime }\pi (a)u^{\prime \ast }=\left[ 
\begin{array}{cc}
\alpha (a) & \beta (a)u^{\ast } \\ 
u\gamma (a) & \alpha (a)%
\end{array}%
\right] $. To ease notations, we set $\beta ^{\prime }:a\in A\mapsto \beta
(a)u^{\ast }$ and $\gamma ^{\prime }:a\in A\mapsto u\gamma (a)$. We also
denote $\mathcal{H}_{1}$ by $\mathcal{J}$ and (up to a trivial isomorphism)
we write $\mathcal{H}=\mathcal{J}\oplus \mathcal{J}$. Now $\alpha ,\beta
^{\prime }$ and $\gamma ^{\prime }$ are all three linear maps on $\mathcal{J}
$. The representation $u^{\prime }\pi u^{\prime \ast }$ is denoted by $%
\theta $.

Let $b\in A_{1}$ and $a\in A_{-1}$. Then $\left( ba\right) ^{2}\in A_{1}$
and:%
\begin{eqnarray*}
&&\left[ 
\begin{array}{cc}
\alpha \left( \left( ba\right) ^{2}\right) & 0 \\ 
0 & \alpha \left( \left( ba\right) ^{2}\right)%
\end{array}%
\right] =\theta \left( \left( ba\right) ^{2}\right) =\left( \theta (b)\theta
(a)\right) ^{2} \\
&=&\left( \left[ 
\begin{array}{cc}
\alpha (b) & 0 \\ 
0 & \alpha (b)%
\end{array}%
\right] \left[ 
\begin{array}{cc}
0 & \beta ^{\prime }(a) \\ 
\gamma ^{\prime }(a) & 0%
\end{array}%
\right] \right) ^{2} \\
&=&\left[ 
\begin{array}{cc}
\alpha (b)\beta ^{\prime }(a)\alpha (b)\gamma ^{\prime }(a) & 0 \\ 
0 & \alpha (b)\gamma ^{\prime }(a)\alpha (b)\beta ^{\prime }(a)%
\end{array}%
\right]
\end{eqnarray*}%
and thus for all $a\in A_{-1}$ and $b\in A_{1}$ we have: 
\begin{equation*}
\alpha (b)\beta ^{\prime }(a)\alpha (b)\gamma ^{\prime }(a)=\alpha (b)\gamma
^{\prime }(a)\alpha (b)\beta ^{\prime }(a)\text{.}
\end{equation*}%
Now, since $\alpha (A_{1})$ is $\limfunc{SOT}$-dense in $\mathcal{B}(%
\mathcal{J})$ we conclude that for all $T\in \mathcal{B}(\mathcal{J})$ we
have for all $a\in A_{-1}$:%
\begin{equation*}
T\beta ^{\prime }(a)T\gamma ^{\prime }(a)=T\gamma ^{\prime }(a)T\beta
^{\prime }(a)
\end{equation*}%
and thus we have $\beta ^{\prime }(a)T\gamma ^{\prime }(a)=\gamma ^{\prime
}(a)T\beta ^{\prime }(a)$ for all $T\in \mathcal{B}(\mathcal{J})$ and $a\in
A_{-1}$. By Lemma (\ref{central}), for each $a\in A_{-1}$ there exists $%
\lambda (a)\in \mathbb{C}$ such that $\lambda (a)\beta ^{\prime }(a)=\gamma
^{\prime }(a)$. On the other hand, let $a,b\in A_{-1}$ be given. Then:%
\begin{eqnarray*}
\lambda (a)\beta ^{\prime }(a)+\lambda (b)\beta ^{\prime }(b) &=&\gamma
^{\prime }(a)+\gamma ^{\prime }(b)=\lambda (a+b)\beta ^{\prime }(a+b) \\
&=&\lambda (a+b)\beta ^{\prime }(a)+\lambda (a+b)\beta ^{\prime }(b)\text{.}
\end{eqnarray*}%
If $\beta ^{\prime }(a)$ and $\beta ^{\prime }(b)$ are linearly independent
then $\lambda (a)=\lambda (b)=\lambda (a+b)$ (thus $\lambda $ is constant if 
$\beta ^{\prime }(A_{-1})$ is at least two dimensional).

If instead, $\beta ^{\prime }(a)=t\beta ^{\prime }(b)$ for some $t\in 
\mathbb{C}$ then we get: 
\begin{equation*}
\lambda (ta)\beta ^{\prime }(ta)=\gamma ^{\prime }(ta)=t\gamma ^{\prime
}(a)=t\lambda (a)\beta ^{\prime }(a)\text{.}
\end{equation*}%
Hence, if $t\not=0$ and $\beta ^{\prime }(a)\not=0$ then $\lambda
(ta)=\lambda (a)$.

Thus, if $a,b\in A_{-1}$ and $a,b$ are not in $\ker \beta ^{\prime }$ then $%
\lambda (a)=\lambda (b)$ (as $\{a,b\}$ is either linearly independent or
they are dependant but $\beta ^{\prime }(a)$ and $\beta ^{\prime }(b)$ are
not zero). We can make the choice we wish for $\lambda (a)$ when $a\in \ker
\beta ^{\prime }$, so naturally we set $\lambda (a)=\lambda (b)$ for any $%
b\in A_{-1}\backslash \ker \beta ^{\prime }$ (note that $A_{-1}\backslash
\ker \beta ^{\prime }\not=\emptyset $ since $\theta $ is irreducible and
since $\beta ^{\prime }(a)=\gamma ^{\prime }(a^{\ast })^{\ast }$ for all $%
a\in A$). With this choice, we have shown that there exists a $\lambda \in 
\mathbb{C}$ such that $\lambda \beta ^{\prime }(a)=\gamma ^{\prime }(a)$ for
all $a\in A_{-1}$.

Moreover, let $a\in A_{-1}$. Then $\beta ^{\prime }(a^{\ast })=\gamma
^{\prime }(a)^{\ast }$ and $\beta ^{\prime }(a)^{\ast }=\gamma ^{\prime
}(a^{\ast })$ by definition of $\beta ^{\prime }$ and $\gamma ^{\prime }$,
yet $\gamma ^{\prime }(a)=\lambda \beta ^{\prime }(a)$. So if $a=a^{\ast }$
then: 
\begin{eqnarray*}
\gamma ^{\prime }(a) &=&\lambda \beta ^{\prime }(a)=\lambda \beta ^{\prime
}(a^{\ast })=\lambda \gamma ^{\prime }(a)^{\ast }=\lambda \left( \lambda
\beta ^{\prime }(a)\right) ^{\ast } \\
&=&\left\vert \lambda \right\vert ^{2}\beta ^{\prime }(a)^{\ast }=\left\vert
\lambda \right\vert ^{2}\gamma ^{\prime }(a)\text{.}
\end{eqnarray*}%
Now, suppose that $\gamma ^{\prime }(a)=0$ for all $a=a^{\ast }\in A_{-1}$.
By assumption, $\gamma ^{\prime }$ is not zero (since then $\beta ^{\prime }$
would be since $\beta ^{\prime }(a)=\gamma ^{\prime }(a^{\ast })^{\ast }$
and then $\theta $ would be reducible), so there exists $a\in A_{-1}$ such
that $a^{\ast }=-a$ and $\gamma ^{\prime }(a)\not=0$ (since $\gamma ^{\prime
}$ linear and every element in $A_{-1}$ is of the sum of a self-adjoint and
anti-selfadjoint element in $A_{-1}$). But then $ia$ is self-adjoint, and
since $\gamma ^{\prime }$ is linear, $\gamma ^{\prime }(ia)=0$. This is a
contradiction. Hence there exists $a\in A_{-1}$ such that $a=a^{\ast }$ and $%
\gamma ^{\prime }(a)\not=0$. Therefore, $\left\vert \lambda \right\vert
^{2}=1$. Let $\eta $ be any square root of $\lambda $ in $\mathbb{C}$.

Set $\nu =\left[ 
\begin{array}{cc}
1 & 0 \\ 
0 & \eta%
\end{array}%
\right] $ and $\psi =\nu \theta \nu ^{\ast }$ so that 
\begin{equation*}
\psi (a)=\left[ 
\begin{array}{cc}
\alpha (a) & \eta \beta ^{\prime }(a) \\ 
\eta \beta ^{\prime }(a) & \alpha (a)%
\end{array}%
\right] \text{.}
\end{equation*}

Let $v^{\prime }=\frac{1}{\sqrt[2]{2}}\left[ 
\begin{array}{cc}
1 & 1 \\ 
1 & -1%
\end{array}%
\right] $ so that: 
\begin{equation*}
v^{\prime }\psi (a)v^{\prime \ast }=\left[ 
\begin{array}{cc}
\alpha (a)+\eta \beta ^{\prime }(a) & 0 \\ 
0 & \alpha (a)-\eta \beta ^{\prime }(a)%
\end{array}%
\right] \text{.}
\end{equation*}%
Letting $\varphi =\alpha +\eta \beta ^{\prime }$ we see that $\varphi $ is a
*-representation of $A$ and that $\pi $ is unitarily equivalent to the
representation $\pi _{\varphi }$ defined by $\pi _{\varphi }(a)=\left[ 
\begin{array}{cc}
\varphi (a) & 0 \\ 
0 & \varphi \left( \sigma (a)\right)%
\end{array}%
\right] $ and $\pi _{\varphi }(W)=\left[ 
\begin{array}{cc}
0 & 1 \\ 
1 & 0%
\end{array}%
\right] $. In particular, $\pi ^{\prime }=\varphi \oplus \varphi \circ
\sigma $ is a reducible representation of $A$.

Note that we could have done the same proof by limiting ourselves to the
case where $a\in A_{-1}$ is selfadjoint and by calculating $\pi (a)^{\ast
}\pi (a)$, using Lemma (\ref{central2}) instead of Lemma (\ref{central}).
\end{proof}

We read from the proof of Theorem (\ref{Rep}) the following description of
some irreducible representations of $A\rtimes _{\sigma }\mathbb{Z}_{2}$
which completes the statement of Proposition (\ref{Rep1}):

\begin{corollary}
Let $\pi $ be an irreducible representation of $A\rtimes _{\sigma }\mathbb{Z}%
_{2}$ and $\pi ^{\prime \prime }$ its restriction to $A_{1}$. Let $\alpha
,\delta $ be the irreducible representations of $A_{1}$ such that $\pi
^{\prime \prime }=\alpha \oplus \delta $. Then the following statements are
equivalent:

\begin{itemize}
\item $\alpha $ is unitarily equivalent to $\delta $,

\item $\pi $ is unitarily equivalent to a representation $\rho $ such that $%
\rho (W)=\left[ 
\begin{array}{cc}
0 & 1 \\ 
1 & 0%
\end{array}%
\right] $ and $\rho (a)=\left[ 
\begin{array}{cc}
\rho ^{\prime }(a) & 0 \\ 
0 & \rho ^{\prime }\circ \sigma (a)%
\end{array}%
\right] $ where $\rho ^{\prime }$ is an irreducible representation of $A$
and $W$ is the canonical unitary in $A\rtimes _{\sigma }\mathbb{Z}_{2}$ and $%
\rho ^{\prime }$ is not unitarily equivalent to $\rho ^{\prime }\circ \sigma 
$.
\end{itemize}
\end{corollary}

We easily observe that both types of representations described in
Proposition (\ref{Rep0}) and Theorem (\ref{Rep}) do actually occur.

\begin{example}
Let $A=M_{2}$ and $\sigma :a\in M_{2}\mapsto WaW$ where $W=\left[ 
\begin{array}{cc}
1 & 0 \\ 
0 & -1%
\end{array}%
\right] $. All irreducible representations of $M_{2}\rtimes _{\sigma }%
\mathbb{Z}_{2}$ are unitarily equivalent to the identity representation of $%
M_{2}$.
\end{example}

\begin{example}
\label{Ex1}Let $A=C(\mathbb{T})$ and $\sigma :f\mapsto f\circ \sigma ^{\ast
} $ where $\sigma ^{\ast }:\omega \in \mathbb{T}\mapsto -\omega $. Then all
irreducible representations of $C(\mathbb{T})\rtimes _{\sigma }\mathbb{Z}%
_{2} $ are given by the construction of Proposition (\ref{Rep0}). Indeed, if 
$\pi ^{\prime }$ is the restriction of an irreducible representation $\pi $
of $C(\mathbb{T})\rtimes _{\sigma }\mathbb{Z}_{2}$ then $\pi ^{\prime }$ is
irreducible if and only if $\pi ^{\prime }$ is one-dimensional. In this
case, $\pi $ is one-dimensional and thus corresponds to a fixed point in $%
\mathbb{T}$ for $\sigma $. Since there is no such fixed point, $\pi ^{\prime
}$ is reducible and the direct sum of the evaluations at $\omega $ and $%
-\omega $ for some $\omega \in \mathbb{T}$.
\end{example}

\begin{example}
Both types of representations occur if we replace $\sigma ^{\ast }$ in
Example (\ref{Ex1}) by $\sigma ^{\ast \ast }:\omega \in \mathbb{T}\mapsto 
\overline{\omega }$. With the notations of Example (\ref{Ex1}), $\pi
^{\prime }$ is irreducible if and only if it is the evaluation at one of the
fixed points $1$ or $-1$. In this case, $\pi (W)=\pm 1$. Otherwise, $\pi
^{\prime }$ is reducible and the direct sum (up to unitary conjugation) of
the evaluations at $\omega $ and $\overline{\omega }$ for $\omega \in 
\mathbb{T}\backslash \left\{ -1,1\right\} $.
\end{example}

We can deduce one more interesting piece of information on the structure of
irreducible representations of $A\rtimes _{\sigma }\mathbb{Z}_{2}$ from the
proof of Theorem (\ref{Rep}):

\begin{corollary}
Let $\pi $ be an irreducible representation of $A$. Then there exists a
unitary $u$ such that $u^{2}=1$ and $u\pi u^{\ast }=\pi \circ \sigma $ if
and only if the restriction $\pi ^{\prime \prime }$ of $\pi $ to the fixed
point C*-algebra $A_{1}$ is the sum of two unitarily non-equivalent
(irreducible) representations of $A_{1}$.
\end{corollary}

\subsection{Representation theory of $A\rtimes _{\protect\sigma }\mathbb{Z}$
with $\protect\sigma ^{2}=\limfunc{Id}$}

We wish to point out that the previous description of the representation
theory of the crossed-product $A\rtimes _{\sigma }\mathbb{Z}_{2}$ can be
used to derive just as well the representation theory of $A\rtimes _{\sigma }%
\mathbb{Z}$, as described in the following proposition. The
C*-crossed-product $A\rtimes _{\sigma }\mathbb{Z}$ is the universal
C*-algebra generated by $A$ and a unitary $W_{\mathbb{Z}}$ with the
relations:\ $W_{\mathbb{Z}}aW_{\mathbb{Z}}^{\ast }=\sigma (a)$ for all $a\in
A$ \cite{Zeller-Meier68}.

\begin{proposition}
Let $\sigma $ be an order-two *-automorphism of a unital C*-algebra $A$. Let 
$\pi _{2}$ be an irreducible representation of $A\rtimes _{\sigma }\mathbb{Z}%
_{2}$ on some Hilbert space $\mathcal{H}$. Let $\lambda \in \mathbb{T}$.
Denote by $W$ the canonical unitary in $A\rtimes _{\sigma }\mathbb{Z}_{2}$
and $W_{\mathbb{Z}}$ the canonical unitary in $A\rtimes _{\sigma }\mathbb{Z}$%
. Set $\pi $ on $A$ by $\pi (a)=\pi _{2}(a)$ for all $a\in A$ and set $\pi
(W_{\mathbb{Z}})=\lambda \pi _{2}(W)$. Then $\pi $ extends uniquely to a
representation of $A\rtimes _{\sigma }\mathbb{Z}$. Moreover, all irreducible
representations of $A\rtimes _{\sigma }\mathbb{Z}$ are obtained this way.
\end{proposition}

\begin{proof}
It is obvious that $\pi $ thus constructed from $\pi _{2}$ is an irreducible
representation of $A\rtimes _{\sigma }\mathbb{Z}$. Let now $\pi $ be an
irreducible representation of $A\rtimes _{\sigma }\mathbb{Z}$. Since $\pi $
is irreducible and $\pi (W_{\mathbb{Z}})^{2}$ commutes with $\pi (A)$ (since 
$\sigma ^{2}=1$), we conclude that $\pi (W_{\mathbb{Z}})^{2}=\lambda ^{2}$
for some $\lambda \in \mathbb{T}$. Let $U=\lambda ^{-1}\pi (W_{\mathbb{Z}})$%
. Then $U$ is an order-two unitary. Define $\pi _{2}(a)=\pi (a)$ for all $%
a\in A$ and $\pi _{2}(W)=U$: by universality of $A\rtimes _{\sigma }\mathbb{Z%
}_{2}$, the map $\pi _{2}$ extends to a representation of $A\rtimes _{\sigma
}\mathbb{Z}_{2}$. It is irreducible since $\pi $ is. This proves our
proposition.
\end{proof}

\section{Application to C*-crossed-products of $C^{\ast }(\mathbb{F}_{2})$}

This section concerns itself with two examples of an action on the free
group $\mathbb{F}_{2}$ on two generators. This paper deals with
representation theory, so we present here examples which can be handled
using representation theory more or less directly. More precisely, given the
universal C*-algebra $C^{\ast }\left( \mathbb{F}_{2}\right) $ generated by
two unitaries $U$ and $V$, there are three obvious and natural automorphisms
of order 2 to consider: $\alpha $ defined by $\alpha (U)=U^{\ast }$ and $%
\alpha (V)=V^{\ast }$, as well as $\beta $ defined by $\beta (U)=-U$ and $%
\beta (V)=-V$ and at last $\gamma $ defined by $\gamma (U)=V$ and $\gamma
(V)=U$. A companion paper \cite{CLat06} to this one by the same authors
deals with the interesting structure of the fixed point C*-algebra for $%
\gamma $, and thus the study of the related C*-crossed-product of $C^{\ast
}\left( \mathbb{F}_{2}\right) $ by $\gamma $ is done in \cite{CLat06} as
well. The study of $\alpha $ and $\beta $ is undertaken in this section.

The following propositions will help us compute the K-theory of these
crossed-products by bringing the problem back to simple type I
crossed-products on Abelian C*-algebras, to which it will be easy to apply
Theorem (\ref{Rep}).

\begin{proposition}
\label{freecross}Let $A_{1}$ and $A_{2}$ be two unital C*-algebras, and let $%
\alpha _{1}$ and $\alpha _{2}$ be two actions of a discrete group $G$ on $%
A_{1}$ and $A_{2}$ respectively. Let $\alpha $ be the unique action of $G$
on $A_{1}\ast _{\mathbb{C}}A_{2}$ extending $\alpha _{1}$ and $\alpha _{2}$.
Then: 
\begin{equation*}
\left( A_{1}\ast _{\mathbb{C}}A_{2}\right) \rtimes _{\alpha }G=\left(
A_{1}\rtimes _{\alpha _{1}}G\right) \ast _{C^{\ast }(G)}\left( A_{2}\rtimes
_{\alpha _{2}}G\right)
\end{equation*}%
where the free product is amalgated over the natural copies of $C^{\ast }(G)$
in $A_{1}\rtimes _{\alpha _{1}}G$ and $A_{2}\rtimes _{\alpha _{2}}G$
respectively.
\end{proposition}

\begin{proof}
This result follows from universality. Since $G$ is discrete, there is a
natural embedding $i_{k}:C^{\ast }(G)\longrightarrow A_{k}\rtimes _{\alpha
_{k}}G$ for $k=1,2$. Now, given a commuting diagram:%
\begin{equation}
\begin{array}{ccc}
C^{\ast }(G) & \overset{i_{1}}{\longrightarrow } & A_{1}\rtimes _{\alpha
_{1}}G \\ 
\downarrow i_{2} &  & \downarrow j_{1} \\ 
A_{2}\rtimes _{\alpha _{2}}G & \overset{j_{2}}{\longrightarrow } & B%
\end{array}
\label{freecomm}
\end{equation}%
by universality of the amalgated free product, there exists a unique
surjection $\varphi _{B}:\left( A_{1}\rtimes _{\alpha _{1}}G\right) \ast
_{C^{\ast }(G)}\left( A_{2}\rtimes _{\alpha _{2}}G\right) \longrightarrow B$
such that, if we use the notations:%
\begin{equation*}
\begin{array}{ccc}
C^{\ast }(G) & \overset{i_{1}}{\longrightarrow } & A_{1}\rtimes _{\alpha
_{1}}G \\ 
\downarrow i_{2} &  & \downarrow \varphi _{1} \\ 
A_{2}\rtimes _{\alpha _{2}}G & \overset{\varphi _{2}}{\longrightarrow } & 
\left( A_{1}\rtimes _{\alpha _{1}}G\right) \ast _{C^{\ast }(G)}\left(
A_{2}\rtimes _{\alpha _{2}}G\right)%
\end{array}%
\end{equation*}

then: $\varphi _{B}\circ \varphi _{k}=j_{k}$ for $k=1,2$. Of course, up to a
*-isomorphism, there is a unique such universal object. Let us prove that $%
\left( A_{1}\ast _{\mathbb{C}}A_{2}\right) \rtimes _{\alpha }G$ is this
universal object, which will prove the proposition.

First, let $g\in G$ and let $U^{g}\in C^{\ast }(G)$, $U_{1}^{g}=i_{1}\left(
U^{g}\right) \in A_{1}\rtimes _{\alpha _{1}}G$ and $U_{2}^{g}=i_{2}\left(
U^{g}\right) \in A_{2}\rtimes _{\alpha _{2}}G$ and $U_{3}^{g}\in \left(
A_{1}\ast _{\mathbb{C}}A_{2}\right) \rtimes _{\alpha }G$ be the naturally
associated unitaries. Now, we observe that $\left( A_{1}\ast _{\mathbb{C}%
}A_{2}\right) \rtimes _{\alpha }G$ fits in the commutative diagram:%
\begin{equation}
\begin{array}{ccc}
C^{\ast }(G) & \overset{i_{1}}{\longrightarrow } & A_{1}\rtimes _{\alpha
_{1}}G \\ 
\downarrow i_{2} &  & \downarrow \theta _{1} \\ 
A_{2}\rtimes _{\alpha _{2}}G & \overset{\theta _{2}}{\longrightarrow } & 
\left( A_{1}\ast _{\mathbb{C}}A_{2}\right) \rtimes _{\alpha }G%
\end{array}
\label{freecomm2}
\end{equation}%
where $\theta _{k}(a)=a$ and $\theta _{k}(U_{k}^{g})=U_{3}^{g}$ for $a\in
A_{k}$ and $k=1,2$. Indeed, one checks immediately that, for $k=1,2$, the
map $\theta _{k}$ satisfies $\theta _{k}(U_{k}^{g})\theta _{k}(a)\theta
_{k}(U_{k}^{g})^{\ast }=\alpha _{k}(a)=\theta _{k}\left( \alpha
_{k}(a)\right) $ and then we can extend $\theta _{k}$ by universality of $%
A\rtimes _{\alpha _{k}}G$. The commutativity of the diagram is obvious.

Now, let us be given a C*-algebra $B$ fitting in the commutative diagram (%
\ref{freecomm}). Let $a\in A_{k}$ ($k=1,2$). Then set $\psi (a)=j_{k}(a)$.
Note that $\psi (1)=j_{1}(1)=j_{2}(1)=j_{k}\circ i_{k}(1)$ as $i_{k}$ is
unital for $k=1,2$. Hence, $\psi $ extends to $A_{1}\ast _{\mathbb{C}}A_{2}$
by universality of $A_{1}\ast _{\mathbb{C}}A_{2}$. Now, with the notations
of (\ref{freecomm2}), we have $\theta _{1}(U_{1}^{g})=\theta _{2}\left(
U_{2}^{g}\right) =U_{3}^{g}$ by construction. We set $\psi
(U_{3}^{g})=j_{1}\left( U_{1}^{g}\right) =j_{1}\circ i_{1}(U^{g})$. As the
diagram (\ref{freecomm}) is commutative, $\psi (U_{3}^{g})=j_{1}\circ
i_{2}(U^{g})$. Moreover, $\psi (U_{3}^{g})\psi (a)\psi (U_{3}^{g})^{\ast
}=j_{k}\left( U_{k}^{g}aU_{k}^{g\ast }\right) =j_{k}\left( \alpha
_{k}(a)\right) $ for all $a\in A_{k}$ with $k=1,2$ by construction of $\psi $%
. It is easy to deduce that $\psi (U_{3}^{g})\psi (a)\psi (U_{3}^{g})^{\ast
}=\psi \circ \alpha (a)$ for all $a\in A_{1}\ast _{\mathbb{C}}A_{2}$. Hence,
by universality of the crossed-product, the map $\psi $ extends to $\left(
A_{1}\ast _{\mathbb{C}}A_{2}\right) \rtimes _{\alpha }G$ into $B$. Moreover,
by construction $\psi \circ \theta _{1}=j_{1}$ and $\psi \circ \theta
_{2}=j_{2}$. Thus, $\left( A_{1}\ast _{\mathbb{C}}A_{2}\right) \rtimes
_{\alpha }G$ is universal for the diagram (\ref{freecomm}), so $\left(
A_{1}\ast _{\mathbb{C}}A_{2}\right) \rtimes _{\alpha }G=\left( A_{1}\rtimes
_{\alpha _{1}}G\right) \ast _{C^{\ast }(G)}\left( A_{2}\rtimes _{\alpha
_{2}}G\right) $.
\end{proof}

\begin{proposition}
\label{KFree}Let $A_{1}$ and $A_{2}$ be two unital C*-algebras with two
respective one-dimensional representations $\varepsilon _{1}$ and $%
\varepsilon _{2}$. Let $\alpha _{1}$ and $\alpha _{2}$ be two actions of a
discrete group $G$ on $A_{1}$ and $A_{2}$ respectively such that $%
\varepsilon _{1}\circ \alpha _{1}=\varepsilon _{1}$ and $\varepsilon
_{2}\circ \alpha _{2}=\varepsilon _{2}$. Let $\alpha $ be the unique action
of $G$ on $A_{1}\ast _{\mathbb{C}}A_{2}$ extending $\alpha _{1}$ and $\alpha
_{2}$. Let $i_{k}$ be the natural injection of $C^{\ast }(G)$ into $%
A_{k}\rtimes _{\alpha _{k}}G$ for $k=1,2$. Then $K_{\ast }\left( \left(
A_{1}\ast _{\mathbb{C}}A_{2}\right) \rtimes _{\alpha }G\right) $ equals to:%
\begin{equation*}
\left( K_{\ast }\left( A_{1}\rtimes _{\alpha _{1}}G\right) \oplus K_{\ast
}\left( A_{2}\rtimes _{\alpha _{2}}G\right) \right) /\ker \left( i_{1}^{\ast
}\oplus (-i_{2}^{\ast })\right)
\end{equation*}%
where for any *-morphism $\varphi :A\longrightarrow B$ between two
C*-algebras $A$ and $B$ we denote by $K_{\varepsilon }(\varphi )$ the lift
of $\varphi $ to the K-groups by functoriality (where $\varepsilon \in
\{0,1\}$).
\end{proposition}

\begin{proof}
Let $k\in \{1,2\}$. Denote by $V_{k}^{g}$ the canonical unitary in $%
A_{k}\rtimes _{\alpha _{k}}G$ for $g\in G$ such that $V_{k}^{g}a\left(
V_{k}^{g}\right) ^{\ast }=\alpha _{k}(a)$ for all $a\in A_{k}$. Identify $%
\varepsilon _{k}(a)$ with $\varepsilon _{k}(a)1\in C^{\ast }(G)$ for all $%
a\in A_{k}$. Then by universality of the crossed-product $A_{k}\rtimes
_{\alpha _{k}}G$ and since $U_{g}\varepsilon _{k}(a)U_{g}^{\ast
}=\varepsilon _{k}(a)=\varepsilon _{k}\circ \alpha _{k}(a)$ for all $a\in
A_{k}$ (the latter equality is by hypothesis on $\varepsilon _{k}$), the map 
$\varepsilon _{k}$ extends to $A_{k}\rtimes _{\alpha _{k}}G$ uniquely with $%
\varepsilon _{k}(V_{k}^{g})=U^{g}$ for all $g\in G$ where $U_{g}$ is the
canonical unitary associated to $g\in G$ in $C^{\ast }(G)$. Note that $%
\varepsilon _{k}$ thus extended is valued in $C^{\ast }(G)$.

Using the retractions $\varepsilon _{k}:A_{k}\rtimes _{\alpha
_{k}}G\longrightarrow C^{\ast }(G)$ for $k\in \{1,2\}$, we can apply \cite%
{Cuntz82} and thus the sequence:%
\begin{eqnarray*}
0 &\longrightarrow &K_{\ast }(C^{\ast }(G))\overset{K_{\ast }(i_{1})\oplus
K_{\ast }(i_{2})}{\longrightarrow }K_{\ast }\left( A_{1}\rtimes _{\alpha
_{1}}G\right) \oplus K_{\ast }\left( A_{2}\rtimes _{\alpha _{2}}G\right) \\
&\longrightarrow &K_{\ast }\left( \left( A_{1}\rtimes _{\alpha _{1}}G\right)
\ast _{C^{\ast }(G)}\left( A_{2}\rtimes _{\alpha _{2}}G\right) \right)
\longrightarrow 0
\end{eqnarray*}%
is exact. This calculates the $K$-groups of $\left( A_{1}\rtimes _{\alpha
_{1}}G\right) \ast _{C^{\ast }(G)}\left( A_{2}\rtimes _{\alpha _{2}}G\right) 
$ which, by Proposition (\ref{freecross}) is the crossed-product $\left(
A_{1}\ast _{\mathbb{C}}A_{2}\right) \rtimes _{\alpha }G$.
\end{proof}

\begin{remark}
In particular, if $A_{1}$ is Abelian then the existence of $\varepsilon _{1}$
is equivalent to the existence of a fixed point for the action $\alpha _{1}$.
\end{remark}

Of course, $C^{\ast }\left( \mathbb{Z}_{2}\right) =C(\mathbb{Z}_{2})=\mathbb{%
C}^{2}$. Thus $K_{1}\left( C^{\ast }\left( \mathbb{Z}_{2}\right) \right) =0$
while $K_{0}\left( C^{\ast }\left( \mathbb{Z}_{2}\right) \right) =\mathbb{Z}%
^{2}$ is generated by the spectral projection of the universal unitary $W$
such that $W^{2}=1$.

Now, we use Proposition (\ref{KFree}) to compute the $K$-theory of two
examples. The key in each case is to explicitly calculate the type I
crossed-products $C(\mathbb{T})\rtimes \mathbb{Z}_{2}$. We propose to do so
using Theorem (\ref{Rep}).

\begin{proposition}
Let $\beta $ be the *-automorphism of $C^{\ast }\left( \mathbb{F}_{2}\right) 
$ defined by $\beta (U)=-U$ and $\beta (V)=-V$. Then:%
\begin{equation*}
K_{0}\left( C^{\ast }(\mathbb{F}_{2})\rtimes _{\beta }\mathbb{Z}_{2}\right) =%
\mathbb{Z}\text{ and }K_{1}\left( C^{\ast }(\mathbb{F}_{2})\rtimes _{\beta }%
\mathbb{Z}_{2}\right) =\mathbb{Z}^{2}\text{.}
\end{equation*}
\end{proposition}

\begin{proof}
Let $z$ be the map $\omega \in \mathbb{T}\mapsto \omega $. Write $\beta
=\beta _{1}\ast \beta _{1}$ where $\beta _{1}(z)=-z$. The crossed-product $C(%
\mathbb{T})\rtimes _{\beta _{1}}\mathbb{Z}_{2}$ is $\mathcal{C}=C(\mathbb{T}%
,M_{2})=C(\mathbb{T})\otimes M_{2}$. Indeed, if we set $\psi (f)(x)\mapsto %
\left[ 
\begin{array}{cc}
f(x) & 0 \\ 
0 & f(-x)%
\end{array}%
\right] $ and $\psi (W)=\left[ 
\begin{array}{cc}
0 & 1 \\ 
1 & 0%
\end{array}%
\right] $ then $\psi $ extends naturally to a *-morphism from $C(\mathbb{T}%
)\rtimes _{\beta _{1}}\mathbb{Z}_{2}$ into $\mathcal{C}$. Moreover, the
range of $\psi $ is the C*-algebra spanned by $\psi (z)$ and $\psi (W)$
which is easily checked to be $\mathcal{C}$ by the Stone-Weierstrass
theorem, so $\psi $ is surjective. It is injective as well: let $a\in \ker
\psi $. If $\pi $ is an irreducible *-representation of $C(\mathbb{T}%
)\rtimes _{\beta _{1}}\mathbb{Z}_{2}$ then by Theorem (\ref{Rep}), $\pi $ is
(up to unitary equivalence) acting on $M_{2}$ by $\pi (f)=\left[ 
\begin{array}{cc}
f(x) & 0 \\ 
0 & f(-x)%
\end{array}%
\right] $ and $\pi (W)=\left[ 
\begin{array}{cc}
0 & 1 \\ 
1 & 0%
\end{array}%
\right] $, for some fixed $x\in \mathbb{T}$. Thus, if $\rho _{x}$ is the
evaluation at $x$ in $\mathcal{C}$ then $\rho \circ \psi =\pi $ and thus $%
\pi (a)=0$. Thus $a=0$ as $\pi $ arbitrary and thus $\psi $ is a
*-isomorphism.

Of course, $K_{\ast }\left( M_{2}(C(\mathbb{T))}\right) =K_{\ast }\left( 
\mathbb{T}\right) $ so $K_{0}\left( C(\mathbb{T})\rtimes _{\beta _{1}}%
\mathbb{Z}_{2}\right) =\mathbb{Z}$ and $K_{1}\left( C(\mathbb{T})\rtimes
_{\beta _{1}}\mathbb{Z}_{2}\right) =\mathbb{Z}$. Moreover, $K_{1}$ is
generated by $z$ while $K_{0}$ is simply generated by the identity of $C(%
\mathbb{T})$. The map $i_{k}:C^{\ast }(\mathbb{Z}_{2})\rightarrow \mathcal{C}
$ maps the generator of $C^{\ast }(\mathbb{Z}_{2})$ to $w$, and thus $%
i_{k}^{\ast }$ maps the two spectral projections of $w$ to $1$. Hence, $%
i_{k}^{0}:\mathbb{Z}^{2}\rightarrow \mathbb{Z}$ is defined by $%
i_{k}(0,1)=i_{k}(1,0)=1$. Thus by Proposition (\ref{freecross}), we have $%
K_{0}\left( C^{\ast }(\mathbb{F}_{2})\rtimes \mathbb{Z}_{2}\right) =\mathbb{Z%
}$ and $K_{1}\left( C^{\ast }(\mathbb{F}_{2})\rtimes \mathbb{Z}_{2}\right) =%
\mathbb{Z}^{2}$.
\end{proof}

\begin{proposition}
\label{alpha}Let $\alpha $ be the *-automorphism of $C^{\ast }\left( \mathbb{%
F}_{2}\right) =C^{\ast }(U,V)$ defined by $\alpha (U)=U^{\ast }$ and $\alpha
(V)=V^{\ast }$. Then:%
\begin{equation*}
K_{0}\left( C^{\ast }(\mathbb{F}_{2})\rtimes _{\alpha }\mathbb{Z}_{2}\right)
=\mathbb{Z}^{4}\text{ and }K_{0}\left( C^{\ast }(\mathbb{F}_{2})\rtimes
_{\alpha }\mathbb{Z}_{2}\right) =0\text{.}
\end{equation*}
\end{proposition}

\begin{proof}
Write $\alpha =a_{1}\ast \alpha _{1}$ where $\alpha _{1}(z)=\overline{z}$
where $z$ is the map $\omega \in \mathbb{T}\mapsto \omega $. Now, the
crossed-product $C(\mathbb{T})\rtimes _{\alpha _{1}}\mathbb{Z}_{2}$ is the
C*-algebra $\mathcal{B}=\left\{ h\in C([-1,1],M_{2}):h(1),h(-1)\text{
diagonal}\right\} $. Indeed, define $\psi (f)\left( t\right) $ for all $f\in
C(\mathbb{T})$ and $t\in \lbrack 0,1]$ by:%
\begin{equation*}
\frac{1}{2}\left[ 
\begin{array}{cc}
f\left( t,-y\right) +f\left( t,y\right) & f\left( t,-y\right) -f\left(
t,y\right) \\ 
f\left( t,-y\right) -f\left( t,y\right) & f\left( t,-y\right) +f\left(
t,y\right)%
\end{array}%
\right]
\end{equation*}%
where $y=\sqrt[2]{1-t^{2}}$. Set $\psi (W)=w=\left[ 
\begin{array}{cc}
1 & 0 \\ 
0 & -1%
\end{array}%
\right] $. Then $w^{2}=1$ and $w\psi (f)w=\psi (\alpha _{1}(f))$ so $\psi $
extends to a unique *-morphism from $C(\mathbb{T})\rtimes _{\alpha _{1}}%
\mathbb{Z}_{2}$ into $\mathcal{B}$. By the Stone-Weierstrass theorem, one
can check that $\psi $ is indeed onto. Last, let $\pi $ be an irreducible
*-representation of $C(\mathbb{T})\rtimes _{\alpha _{1}}\mathbb{Z}_{2}$. If
the restriction $\pi ^{\prime }$ of $\pi $ to $C(\mathbb{T})$ is
irreducible, then $\pi ^{\prime }$ is one-dimensional and there exists $x\in 
\mathbb{T}$ such that $\pi ^{\prime }(f)=\pi (f)=f(x)$ for all $f\in C(%
\mathbb{T})$. By Theorem (\ref{Rep}) since $\pi ^{\prime }$ is irreducible, $%
\pi $ is also one-dimensional and $\pi (W)$ is a scalar unitary (hence it is 
$1$ or $-1$ since $W^{2}=1$), so it commutes with $\pi (f)$ for all $f$.
Since $\left( WfW\right) (x)=f(\overline{x})$ we conclude that $x=1$ or $%
x=-1 $. Either way let $\rho _{x}$ be the evaluation at $x$ in $\mathcal{B}$%
. Then $\rho _{x}(h)$ is diagonal by definition of $\mathcal{B}$ for all $%
h\in \mathcal{B}$. Let $\rho _{x,+1}$ be the one-dimensional representation
defined by the upper-left corner of $\rho _{x}$ and let $\rho _{x,-1}$ be
the one-dimensional representation defined by the lower-right corner of $%
\rho _{x}$. Note that either way, $\rho _{x,+}(\psi (f))=\rho _{x,-}(\psi
(f))=f(x)$ for all $f\in C(\mathbb{T})$. On the other hand, $\rho
_{x,\varepsilon }\left( \psi (W)\right) =\varepsilon $. Hence, we have
proven that $\rho _{x,W}\circ \psi =\pi $. Thus if $a\in C(\mathbb{T}%
)\rtimes _{\alpha _{1}}\mathbb{Z}_{2}$ and $\psi (a)=0$ then $\pi (a)=0$.

If instead, $\pi $ restricted to $C(\mathbb{T})$ is reducible, then by
Theorem (\ref{Rep}) $\pi $ is unitarily equivalent to a representation $\pi
^{\prime }$ acting on $M_{2}$ defined as follows: there exists $x\in \mathbb{%
T}$ such that $\pi ^{\prime }(f)=\left[ 
\begin{array}{cc}
f(x) & 0 \\ 
0 & f(\overline{x})%
\end{array}%
\right] $ for all $f\in C(\mathbb{T})$ and $\pi ^{\prime }(W)=\left[ 
\begin{array}{cc}
0 & 1 \\ 
1 & 0%
\end{array}%
\right] $. Up to conjugating by the unitary $\frac{1}{\sqrt[2]{2}}\left[ 
\begin{array}{cc}
1 & -1 \\ 
1 & 1%
\end{array}%
\right] $, we see that if we set $\rho (h)=h(t)$ for all $h\in \mathcal{B}$
where $t$ is defined by $x=\left( t,\sqrt[2]{1-t^{2}}\right) $ then $\rho
\circ \psi =\pi ^{\prime }$ and thus, if $a\in \ker \psi $ then $\pi
^{\prime }(a)=0$ so $\pi (a)=0$. In conclusion, if $a\in \ker \psi $ then $%
\pi (a)=0$ for all (irreducible)\ *-representation of $C(\mathbb{T})\rtimes
_{\alpha _{1}}\mathbb{Z}_{2}$ and thus $a=0$, so $\psi $ is a *-isomorphism.

The $K$-theory of $\mathcal{B}$ is easy to calculate. We start with the
exact sequence $0\rightarrow C_{0}((-1,1),M_{2})\overset{i}{\rightarrow }%
\mathcal{B}\overset{q}{\rightarrow }\mathbb{C}^{4}\rightarrow 0$ where $i$
the inclusion and $q$ the quotient map, also defined by $q(a)=a(1)\oplus
a(-1)$ for $a\in \mathcal{B}$ and identifying the diagonal matrices in $%
M_{2} $ with $\mathbb{C}^{2}$. We also used the notation $C_{0}\left(
X\right) $ for the space of continuous functions on a locally compact space $%
X$ vanishing at infinity. The associated six-terms exact sequence is then
simply:%
\begin{equation*}
\begin{array}{ccccc}
K_{0}\left( C(-1,1)\right) =0 & \overset{K_{0}(i)}{\longrightarrow } & 
K_{0}\left( \mathcal{B}\right) & \overset{K_{0}(q)}{\longrightarrow } & 
\mathbb{Z}^{4} \\ 
\uparrow &  &  &  & \downarrow \delta \\ 
0 & \overset{K_{1}(q)}{\longleftarrow } & K_{1}\left( \mathcal{B}\right) & 
\overset{K_{1}(i)}{\longleftarrow } & \mathbb{Z}=K_{1}\left( C_{0}\left(
-1,1\right) \right)%
\end{array}%
\text{.}
\end{equation*}

The generator of the $K_{1}$ group of $C_{0}(-1,1)\otimes M_{2}$ is the
unitary $u_{1}:t\in (-1,1)\mapsto \exp \left( i\pi t\right) 1_{2}$ where $%
1_{2}$ is the unit of $M_{2}$. However, $u_{1}$ is trivial in $K_{1}(%
\mathcal{B})$ via the obvious homotopy $\left( u_{\lambda }\right) _{\lambda
\in \lbrack 0,1]}$ with $u_{\lambda }:t\in \lbrack -1,1]\mapsto \exp \left(
\pi i\lambda t\right) $ (note that $u_{\lambda }$ for $\lambda \in (0,1)$ is
not in the unitalization of $C_{0}\left( -1,1\right) $ since $u_{\lambda
}(-1)\not=u_{\lambda }(1)$). Thus $K_{1}(i)=0$, $K_{1}(\mathcal{B})=0$ and
the range of $\delta $ is $\mathbb{Z}$ by exactness. Hence, again by
exactness, $\ker \delta $ is a copy of $\mathbb{Z}^{3}$ inside of $\mathbb{Z}%
^{4}=K_{0}\left( \mathbb{C}^{4}\right) $.

Let $p=\frac{1}{2}\left( W+1\right) $ and $p^{\prime }=\frac{1}{2}\left(
1+Wz\right) $ (note that $WzWz=\overline{z}z=1$ so $p^{\prime }$ is a
projection). We calculate easily that $K_{0}\left( q\right) (p)=(1,0,1,0)$
while $K_{0}\left( q\right) (p^{\prime })=(1,0,0,1)$.

The subgroup of $\mathbb{Z}^{4}$ generated by $(1,0,1,0)$, $(1,1,1,1)$ and $%
(1,0,0,1)$ is isomorphic to $\mathbb{Z}^{3}$. By exactness, it must be $%
\mathbb{Z}^{3}$. Since $K_{0}(i)=0$, the map $K_{0}(1)$ is an injection and
thus $K_{0}\left( C(\mathbb{T})\rtimes _{\alpha _{1}}\mathbb{Z}_{2}\right) =%
\mathbb{Z}^{3}$ generated by the spectral projections of $w$ and of $Wz$,
and $K_{1}\left( C(\mathbb{T})\rtimes _{\alpha _{1}}\mathbb{Z}_{2}\right) =0$%
.

Moreover, $i_{k}:C^{\ast }(\mathbb{Z}_{2})\rightarrow C(\mathbb{T})\rtimes
_{\alpha _{1}}\mathbb{Z}_{2}$ maps the generator of $C^{\ast }(\mathbb{Z}%
_{2})$ to $w$, so the range of $i_{k}^{\ast }$ is the subgroup generated by $%
[p]$ and $[1]$. Thus, by Proposition\ (\ref{freecross}), $K_{0}(C^{\ast }(%
\mathbb{F}_{2})\rtimes _{\alpha _{1}}\mathbb{Z}_{2})=\mathbb{Z}^{4}$ and $%
K_{1}\left( C^{\ast }(\mathbb{F}_{2})\rtimes _{\alpha _{1}}\mathbb{Z}%
_{2}\right) =0$.
\end{proof}

\bibliographystyle{amsplain}
\bibliography{thesis}

\end{document}